\documentstyle[11pt, amssymb]{article}
\input{epsf}

\newcommand{\ci}{continuous}
\newcommand{\ru}{ unitary representation}
\newcommand{\Qa}{ Quantum 'ax+b' group }
\newcommand{\qa}{ quantum 'ax+b' group }

\newcommand{\slw}{S.L. Woronowicz}

\newcommand{\cH}{{\cal H}}
\newcommand{\K}{{\cal K}}
\newcommand{\cK}{{\cal K}}
\newcommand{\sH}{_{\cal H}}
\newcommand{\sK}{_{\cal K}}

\newcommand{\rf}[1]{{\rm (\ref{#1})}}
\newcommand{\ch}{ {\cal C}(H)}

\newcommand{\fh}{F_{\hbar}} 
\newcommand{\hb}{\hbar}
\newcommand{\vt}{V_{\theta}} 
\newcommand{\za}{-\!\!\circ} 
\newcommand{\R}{{\mathbb R}}
\newcommand{\C}{{\mathbb C}}

\newcommand{\be}{\beta}
\newcommand{\de}{\delta}

\newcommand{\N}{{\mathbb N}}
\newcommand{\M}{{\rm M}}
\newcommand{\B}{{\rm B}}

\newcommand{\fil}{\varphi}

\newcommand{\qed}{$\Box$}
\newcommand{\Ci}{C_{\infty}}
\newcommand{\Cir}{C_{\infty}(\R)}
\newcommand{\eh}{e^{\frac{i\hbar}{2}} }

\newcommand{\Lkr}{L^{2}(\R)}
\newcommand{\po}{\hat{p}}
\newcommand{\qo}{\hat{q}}

\newcommand{\csta}{$C^{*}$-}
\newcommand{\cstal}{$C^{*}$-algebra}
\newcommand{\te}{\otimes}

\newcommand{\id}{{\rm id}}
\newcommand{\Mor}{{\rm Mor}}

\newcommand{\sign}{{\rm sign }\;}
\newcommand{\whe}{\hspace*{5mm}\mbox{\rm where}\hspace{5mm}}

\newcommand{\af}{\hspace*{1mm} {\bf \eta} \hspace{1mm}}
\newcommand{\od}{\hspace*{5mm}}

\newcommand{\mlot}{\mbox{$\hspace{.5mm}\bigcirc\hspace{-3.7mm}
\raisebox{-.7mm}{$\top$}\hspace{1mm}$}}

\newcommand{\dow}{{\bf Proof: }}

\newcommand{\zw}{CB(\Lkr)}

\newcommand{\mi}{\hspace*{3mm} {\rm and} \hspace{3mm}}

\newcommand{\dla}{\hspace*{5mm} {\rm for} \hspace{5mm}}

\newcommand{\dlad}{for any }

\newcommand{\bfa}{\begin{fakt}}\newcommand{\efa}{\end{fakt}}
\newcommand{\ble}{\begin{lem}}\newcommand{\ele}{\end{lem}}
\newcommand{\bst}{\begin{stw}}\newcommand{\est}{\end{stw}}
\newcommand{\bde}{\begin{defi}}\newcommand{\ede}{\end{defi}}
\newcommand{\bwn}{\begin{wn}}\newcommand{\ewn}{\end{wn}}
\newcommand{\buw}{\begin{uwaga}}\newcommand{\euw}{\end{uwaga}}
\newcommand{\bdy}{\begin{dygresja}}\newcommand{\edy}{\end{dygresja}}
\newcommand{\bwa}{\begin{warning}}\newcommand{\ewa}{\end{warning}}
\newcommand{\bpr}{\begin{przy}}\newcommand{\epr}{\end{przy}}
\newcommand{\btw}{\begin{tw}}\newcommand{\etw}{\end{tw}}
\newcommand{\beq}{\begin{equation}}\newcommand{\eeq}{\end{equation}}
\newcommand{\bit}{\begin{itemize}}\newcommand{\eit}{\end{itemize}}
\newcommand{\bq}{\begin{quote}}\newcommand{\eq}{\end{quote}}
\newcommand{\ba}{\begin{array}}\newcommand{\ea}{\end{array}}

\newtheorem{defi}{Definition}[section]
\newtheorem{wn}[defi]{Observation}
\newtheorem{tw}[defi]{Theorem}
\newtheorem{lem}[defi]{Lemma}

\newtheorem{fakt}[defi]{Corollary}
\newtheorem{stw}[defi]{Proposition}
\newtheorem{przy}[defi]{Example}
\newtheorem{uwaga}[defi]{Remark}
\newtheorem{dygresja}[defi]{Dygresja}
\newtheorem{warning}[defi]{Warning}

\title{Unitary representations of the quantum ``ax+b'' group}
\author{Ma{\l }gorzata  Rowicka - Kudlicka
\thanks{Supported by KBN grant 
No 2 PO3A 036 18}\\
Institute of Mathematics, Polish Academy of Sciences,\\
\'Sniadeckich 8, 00-950 Warszawa, Poland\\
email: rowicka@fuw.edu.pl}
\date{January 31, 2001}
\begin{document}   
\maketitle
\begin{abstract} 
We find all unitary representations of the  
 quantum ``ax+b'' group. It turns out that this quantum group is 
selfdual in the sense that all unitary representations are  
 'numbered' by elements of the same group. Moreover, 
 we discover the  formula for all unitary representations 
involving \slw 's quantum exponential function.
{\bf keywords:} \cstal -- crossed product, quantum group\\
{\bf MSC-class:} 20G42 (Primary), 47L65, 81R15 (Secondary).
\end{abstract}

\section{Introduction}

Locally compact quantum groups are nowadays studied extensively 
by many scientists \cite{ks}.Although at the moment there 
is no commonly accepted definition of a locally compact quantum group,
there are promising approaches and interesting examples have been worked out. 
Among the most interesting ones is the quantum ``ax+b''
 group constructed by S.L. Woronowicz and S. Zakrzewski in \cite{ax+b}.
According to the recent computation by A. van Daele, 
this group is the first known example of an interesting
phenomenon  foreseen by 
Vaes and Kustermans in \cite{vaesk}.
In this paper we study the \qa      
  from the point of view of unitary representations and duality theory.

The aim of this paper is to derive a formula for all unitary 
representations of the \qa.
The \cstal \ of all continuouos functions 
vanishing at infinity on this quantum group  
is 'generated' (in the sense we explain below) 
by three unbonded operators $\log a$, $b$ and $i\beta b$.
We discuss the \qa  on the Hilbert space level in Section 3.
For the reader's convienience in Section 3 and 4 we 
recall the definition and relevant information about the  
\qa \  and our previous  results on unitary representations 
 of some braided quantum groups related to the \qa  \cite{paper1} . 
In Section \ref{cp} we discuss    the \csta \ and 
$W^*$- crossed  products  
connected with \qa\ 
 and prove some 
 propositions we will use later on.
In Section \ref{ost} 
 the formula for  all unitary representations of \qa \ 
is proven. 
In our proof we follow the method used by \slw \ in 
the case of quantum $E(2)$ 
group (see \cite{oe2}).
In the next Section we will fix the  notation.

\section{Notation}
\label{not}
Let 
$\qo$ and $\po$  respectively
denote the canonical coordinate and momentum operators acting on 
 the Hilbert space $\Lkr$. 
The domain of  
the operatora $\qo$ is the set 
$$D(\qo)=\{\psi\in \Lkr\;:\; \int_{\R}x^2|\psi(x)|^2dx<\infty\;\}$$
 and  
for any $\psi\in D(\qo)$ the operator $\qo$ is given  by
\[(\qo \psi)(t)=t\psi(t)\ .\] 
The domain of the  operator $\po$ is the set
$$D(\po)=\{\psi\in \Lkr\;:\; \psi'\in\Lkr\;\}\ ,$$ 
 where distributional differentiation is understood.   

The operator $\po$ is given for 
 any  $\psi$ from  $D(\po)$ by 
$$(\po f)(x)=\frac{{\hb}}{i}\frac{d f(x)}{dx}
,$$
 where distributional differentiation is understood. 

We consider only concrete \cstal s, i.e. embedded into  \cstal
\ of all bounded operators acting  on Hilbert space $\cH$, denoted by $B(\cH)$.
The \cstal
\  of all compact  operators acting  on  $\cH$ will be  denoted by $CB(\cH)$.
All algebras we consider are seperable with the exception of multiplier algebras (see definition of multiplier algebra below).

Let $A$ be \cstal. Then $M(A)$ will denote the
{\em multiplier algebra} 
of $A$, i.e.
\[M(A)=\{m\in B(\cH):\od ma,am\in A\  \ \mbox{\dlad} \ a\in A\}\ .\]
Observe that $A$ is an  ideal in $M(A)$. If   $A$ is a unital 
 \cstal\ , then $A=M(A)$, in general case $A\subset M(A)$.
For example the multiplier algebra of $CB(\cH)$ is the algebra  $B(\cH)$ and 
the multiplier algebra of \csta algebra $\Cir$ 
of all continuous vanishing at infinity functions on $\R$ is the algebra of all continuous bounded functions on $\R$ denoted by $C_{\scriptstyle bounded}(\R)$.
The natural topology on $M(A)$ is the strict topology, 
i.e. we say that a sequence $(m_n)_{n\in\N}$ of $m_n\in M(A)$ 
converges strictly to 0 if for every $a\in A$, we have 
$||m_n a||\rightarrow 0$ and $||a m_n ||\rightarrow 0$, when 
$n\rightarrow +\infty$. 
Whenever we will consider continuous maps from or into 
$M(A)$, we will mean this topology. 

For any \cstal s $A$ and $B$, we will say that $\phi$ is a morphism 
and write $\phi\in \Mor (A,B)$ if  $\phi$ is a * - algebra
homomorphism acting from $A$ into $M(B)$ and such that 
$\phi (A)B$ is dense  in $B$.
Any $\phi\in \Mor (A,B)$ admits unique extension to a * - algebra
homomorphism acting from $M(A)$ into $M(B)$. For any $S\in M(A)$,
operator $\phi(S)$ is given by 
\[\phi(S)(\phi(a)b)=\phi(Ta)b\ ,\]
where $a\in A$ and $b\in B$.

For any closed  operator $T$ cting on   $\cH$ 
 we define its  {\em $z$-transform} by
\[z_T=T(I+T^*T)^{-\frac{1}{2}}\ .\]
\label{ztransf}
Observe that   $z_T\in \B(\cH)$ and  $||z_T||\leq 1$.
Moreover, one can recover  $T$ from $z_T$
\[T=z_T(I-z_T^*z_T)^{-\frac{1}{2}}\ .\]

A closed operator $T$ acting on $A$ is {\em affiliated} with a \cstal \ $A$ 
iff $z_T\in \M (A)$ and  $(I-z_T^*z_T)A$ 
is dense in $A$. 
 A set of all elements affiliated with $A$ is denoted by  
$A^{\af}$\index{$A^{\af}\;\;$\}.
If  $A$ is a unital \cstal , then  $A^{\af}=M(A)=A$, 
 in general case 
\[A\subset M(A)=\{T\af A\;:\;||T||<\infty \}\subset A^{\af} \ .\]
The set of all elements affilated with  $\Cir$ is the set of all  
continuos functions on 
real line $C(\R)$, and a 
 set of all  elements affiliated  with \cstal\ $CB (\cH)$ is a set 
 of all closed operators  $\ch$. 
This last example shows that a product and a sum of two
 elements affiliated with  $A$ may not  be affiliated 
 with $A$, since it is well known that a sum and a product of two 
closed operators may not be closed. Affiliation relation in \cstal\  theory 
  was introduced by   Baaj and  Julg in  \cite{baajaf}. 

Observe, that if  $\phi\in \Mor (A,B)$, 
then one can extend  $\phi$ 
 to elements  affiliated with  $A$. Let us start with the  observation, 
 that for any 
$T\in M(A)$ we have 
\[ \phi (z_T)=z_{\phi (T)}\ .\]
 Hence for any  $T\af A$ we have  $z_T\in M(A)$. Moreover, there exists   
  a unique closed  operator S such that  $\phi (z_T)=z_S$. 
This operator is given by  
$$S=\phi (z_T)\phi (I-z_T^*z_T)^{-\frac{1}{2}}.$$
 From now on we will  write  
 $S=\phi (T)$. 

 We recall now a  nonstandard notion of  generation we use in this paper.  
This notion was introduced in  \cite{wunb}, where  a  
 generalisation of the theory of unital  \cstal s generated by  
 a finite number of generators was presented. It was proved in 
\cite{oper} that such \cstal s  
 are isomorphic to  algebras of  
 continuous operator functions on compact  operator domains
(see  Section 1.3 of  \cite{paper1} and references therein). 
In this approach, the  algebra of all  continuous vanishing at infinity 
 functions on a compact quantum    group  is generated  by matrix 
elements of  fundamental representation. 
  To use this approach to noncompact quantum  groups,
  one has to extend  the notion of  a generation of a  \cstal\  to   
 nonunital \cstal s  and  unbounded generators . 
According to the definition we recall below,  
{\em \cstal\ of continuous  vanishing at infinity functions 
 on a locally  compact quantum  group}
 is  generated  by its fundamental representation . 
However in this case the fundamental  representation is not 
  unitary and  the generators are unbounded  operators,  
 so they are not in the \cstal\  $A$.  

Assume for a while, that were are  given  a \cstal\  $A$  and operators 
 $T_1,T_2,...,T_N$ affiliated with   $A$.
We say that $A$ is {\em  generated} by  $T_1,T_2,...,T_N$ 
if for any Hilbert space 
$\cH$,  a nondegenerate \cstal\ $B\subset B(\cH)$ 
and any $\pi\in Mor(A,CB(\cH))$ we have 
\[\left(
\begin{array}{c} 
\pi (T_i) \mbox{ is affiliated with } A\\
\mbox{ for any } i=1,...,N  
\end{array}
\right)
\Longrightarrow
\left(
\begin{array}{c}
\pi\in \Mor(A,B)
\end{array}
\right).\] 

We stress that described above 'generation' is a relation 
between $A$ and some operators $ T_1,T_2,...,T_N$ and 
both have to be known in advance. There is no 
procedure to obtain $A$ knowing only $ T_1,T_2,...,T_N$ and 
it is even possible that  there is no  $A$ generated by such operators.

For unital  \cstal s  generation in the sense introduced above 
 is exactly the same as the classical notion of generation. 
More precisely, 
  let  $A$ be a unital  \cstal   \ and let  
\mbox{$T_1, T_2, \dots , T_N\in A$.}
If  $A$ is the norm closure of all linear  combinations of  
\mbox{$ I, T_1,  \dots , T_N$,}
 then  $A$ is generated  by $ T_1, T_2, \dots , T_N$ in the sense of the  
 above  definition. 
On the other hand, let    $A$ be a  \cstal\ generated  by 
 $T_1, T_2, \dots , T_N\af A$, such that  $||T_i||<\infty $ for 
$i=1,2,\dots, N$. Then $A$ contains unity , 
$T_1, T_2, \dots , T_N\in A$ and  $A$ is the norm closure of the set of all 
  linear  combinations of 
$ I, T_1, T_2, \dots , T_N$.

An easy example of this relation is that \csta algebra $\Cir$ 
of all continuous vanishing at infinity functions on $\R$ 
is generated by function $f(x)=x$ for any $x\in \R$.
The other example is \cstal \ $\zw$ which is generated by $\po$ and $\qo$.

Let   $A$ and $B$ be  \cstal s and assume that we know generators of $A$. 
In order to describe  $\phi \in \Mor(A,B)$ uniquely it is enough 
  to know how  $\phi$ acts on  generators of  $A$.

We will use exclusively the minimal tensor product of \cstal s  and it will 
be denoted by $\te$.
We will also use the leg numbering notation. For example, if $\phi
\in \Mor (A\te A,A\te A)$ then $\phi_{12}(a\te b)= a\te b\te I_A$
and  $\phi_{13}(a\te b)= a\te I_A\te b$ for $a,b\in A$.
Clearly, $\phi_{12}, \phi_{13}\in \Mor (A\te A,A\te A\te A)$.

Let  $f$ and $\phi$ be strongly commuting selfadjoint operators.
Then, by the spectral theorem
\[f=\int _{\Lambda} r dE(r,\rho)\mi \phi=\int _{\Lambda} \rho
 dE(r,\rho)\ ,\]
where $dE(r,\rho)$ denotes the common spectral measure assiociated 
with $f$ and $\phi$ and $\Lambda$ stands for a joint spectrum 
of $f$ and $\phi$.
Then 
\[
F (f,\phi)=\int _{\Lambda} F(r,\rho) dE(r,\rho)\ .
\]

Let $b$ be a selfadjoint operator and  let  the symbol $\chi$ 
denote the characteristic function defined on  $\R$. 
By $\chi(b\neq 0)$ we mean  the projection  operator on 
 the subspace $\ker b^\perp$, by $\chi(b< 0)$ - the projection  onto  
 the subspace on which $b$ is negative, and so on.

\section{\Qa\ on Hilbert space level}

Now we will introduce the commutation  relations  describing the 
quantum "ax+b" group.
One may guess that only two selfadjoint operators 
  $a$ and  $b$ will appear in such relations. 
However, it was found in  \cite{qexp} and \cite{ax+b}, 
 that in order to ensure existence of a selfadjoint extension 
 of a  sum $a+b$ one has to add an additional 
generator, denoted here by  $ib\beta$.

First we will focus on the operators $b$ and $\beta$.
Operator $\beta$ 
 is an analogue of a noncontinuous  (but measurable) 
 function on the quantum   "ax+b" group, 
so its not a good choice for generator, which should 
 be an analogue of a continuous function.
(For the notion of generator used here see 
Section 2.)
 Namely, the operator $\beta$  itself  will not make a good generator, 
since it is not affiliated with the \cstal\ of all 
continuous functions vanishing at infinity on the 
quantum "ax+b" group, which  we denote by $A$ .
Operator  $i b\beta$ is affiliated with this algebra, 
but $\beta$ is only in  the von Neumann algebra
$ A^{\prime\prime}$.
That is the  reason why we use $ib\beta$ instead of $\beta$. 
 We will  use the  {\em Zakrzewski relation} indroduced by S.L. 
Woronowicz in \cite{qexp}.
We  say that two selfadjoint operators $b$ and $d$ 
satisfy the Zakrzewski commutation relation 
and we write 
$$b\za d \ ,$$
 if $b$ commutes with $\sign d$ and 
\beq
|b|^{it}d|b|^{-it}=e^{\hb t}d \dla t\in \R\ ,
\label{zakrel}
\eeq
where $-\pi<\hb<\pi$.
Here $\sign b$ denotes the selfadjoint, bounded operator 
 appearing in the polar decomposition of $b$.
Any pair of operators $b$ and $d$ such that $b\za d$ and 
$\ker b=\ker d=\{0\}$ is 
unitarily equivalent to a direct sum of a certain number of copies of 
$\pm e^{\qo}$ and $\pm e^{\po}$.

Consider  a pair of selfadjoint 
operators $(b,\be)$  
satisfying  
\beq
 b\be=-\be b\od \be^2=\chi(b\neq 0)\ ,
\label{warM}
\eeq 
where 
 the second condition means that $\be^2$ 
 is a projection  operator onto  
 the subspace $\ker b^\perp$ (see Section 2).  
If a pair of selfadjoint operators 
$(b,\be)$  acting on   $\cH$ satisfies \rf{warM} we will write 
$$(b,\be)\in M_\cH\ $$
and call $(b,\be)$ an {\em $M$-pair}.

$M$-pairs were studied in detail in \cite{paper1}.
If two pairs 
$(b,\be)\in M_\cH$ and $(d,\de)\in M_\cH$
 satisfy an additional  condition
\beq
b\za d \od b \de=\de b\od d\be=\be d\od \be\de=\de\be\ ,
\label{warzgo}
\eeq
we write
\[ (b,\be,d,\de) \in M^2 _{\cH} \  .\]

For $(b,\be,d,\de) \in M^2 _{\cH}$ we will  define
the $\mlot$ product of M-pairs 
\beq
\label{mlot}
(b,\be)\mlot (d,\de)= (\tilde{d},\tilde{\delta})\ .
\eeq

Before we write down  formulas for $(\tilde{d},\tilde{\delta})$, we 
will give 
some explanation  why our choice is natural.
First of all, we want $(\tilde{d},\tilde{\delta})$ to be an 
M-pair acting on $\cH$.
Secondly, having already in mind application to the \qa, we  define 
this product in such a way that the first element is as close to $b+d$ 
as possible. The reason why this is so desirable will be obvious in 
Section \ref{ost}.

Here, however, we encounter a serious difficulty.
It is  well known that the sum of  two selfadjoint operators need 
 not to be selfadjoint, so in general the simplest proposal
\[\tilde{d}=b+d\]
will not do. The second guess  is to
 take a selfadjoint extension of the sum above. 
However, in some cases, namely when defficiency  
indices of $b+d$ are not equal, 
there are no selfadjoint extensions. 
An easy  example is to take  $b=e^{\qo}$ and 
$d=-e^{\po}$. Then $b+d$ is a closed symmetric operator, but there 
is no selfadjoint 
extension since deficiency indices are 0 and 1.

The theory of selfadjoint extensions of sums of operators 
satisfying Zakrzewski relation   was 
developed   by S.L. Woronowicz in \cite{qexp}. 
It was proved there that every selfadjoint extension of $b+d$ defines 
uniquely a selfadjoint operator $\phi$ such that $\phi$ anticommutes with $b$ and $\delta$ 
and such that $\phi ^2=\chi(\eh bd <0)$.

Let $\fh$  denote the  quantum exponential function 
indroduced by S.L. Woronowicz in \cite{qexp} as 
\beq
\fh(r,\rho)=\left\{
\begin{array}{ccc}
\vt (\log r)& \mbox{\rm dla } & r>0 \;\;\;\mbox{\rm  i } \rho = 0\\
\{1+i\rho |r|^{\frac{\pi}{\hbar}}\}\vt(\log |r|-\pi i)& 
\mbox{\rm dla} & r<0 \;\;\;\mbox{\rm i } \rho = \pm 1
\end{array}\ ,
\right.
\eeq
where $\theta=\frac{2\pi}{\hb}$ and  $\vt$ is a meromorphic function 
on $\C$ defined as 
\[\vt(x)=\exp\left\{\frac{1}{2\pi i}\int_0^\infty \log(1+a^{-\theta})
\frac{da}{a+e^{-x}}\right\}\ .\]

From now on we assume that 
\beq
\label{hb}
\hb = \pm \frac{\pi}{2k+3}, \whe k=0,1,2,...
\eeq

The reason why we restrict ourselves 
only to such values of parameter $\hb$ is that 
 we are motivated by the \qa and 
 and it turns out that only for such $\hb$ the 
 \qa exists on \cstal\   and Hilbert space level (see \cite{ax+b}).
 Without  loss of generality we may  assume that $\ker d=\{0\}$, 
since in case $d=0$ obviously there are no problems 
with selfadjointness of $b+d=b$.
 If $d$ is invertible then the operator $\eh d^{-1}b$ is 
selfadjoint \cite{qexp}. In such a case   a selfadjoint 
extension of  $b+d$ is given by 
\beq
\label{ext}
 \tilde{d}=[b+d]_{\phi }=F_{\hbar}(f,\phi)^* d
F_{\hbar}(f,\phi)\ ,
\eeq
where 
\[f= \eh d^{-1}b \mi \phi=(-1)^k\be\de\chi( \eh bd<0)\ \]
and the relation between $\hb$ and $k$ is given by \rf{hb}.
Analogously
\[
\tilde{\delta}=F_{\hbar}(f,\phi)^*\delta 
F_{\hbar}(f,\phi)\ .\]

Using Woronowicz's results  we proved in \cite{paper1}  
that if $(b,\be,d,\de)\in M^2 _{\cH}$ 
then  a selfadjoint extension of $b+d$  always exists and may be 
 given by \rf{ext}. 
 Moreover, then  $([b+d]_{\phi},\tilde{\delta})\in M_H$.
We also  proved in that paper (Theoreme 3.3), that  
\ble
\label{lem}
$U$  is a unitary representation of $M$ 
acting on a Hilbert space $\cK$, 
i.e. $U$ is an  operator map (see  Section 1.3 of  
\cite{paper1} and references therein)
$$ U\ :\ M_{\cH}\ \longrightarrow 
B(\cK\te \cH)$$ 
satisfying 
\[ U(b,\be) U (d,\de) =U\left((b,\be)\mlot (d,\de)
\right)\ ,\]
iff there is $(g,\gamma)\in M_{\cal K}$ such that 
for any $(b,\beta)\in M_{\cal H}$
\beq
\label{repM}
U( b,\be)=\fh(g\te b,(\gamma\te \be)
\chi(g\te b<0)).
\eeq
\ele
It means, that unitary representations of $M$ acting on $\cK$ are 
'numbered' by pairs $(g, \gamma)\in M_{\cal K} $.
This Lemma  will be of great use in the proof of  Theorem \ref{ru}, which is 
 the main result of this paper.

Now we are ready to introduce commutation relations related to the 
\qa .
We say that
$(a,b,\be)\in G_{\sH}$
 and we call $(a,b,\be)$ a $G$-{\em triple}  if  
  $a$,\  $b$ and $\beta$  are selfadjoint   operators
 acting  on  a Hilbert space   $\cH$,   $a$ is  positive and 
  invertible , $a\za b$,\   $a$ commutes with 
  $\be$  and $(b,\be)\in M_{\cH}$.
One can define operation $\mlot$ on  $G$-triples 
 for any $(a,b,\be)\in G_{\sH}$ and $(c,d,\de)\in G_{\sK}$
\[(a,b,\be)\mlot (c,d,\de)=(\tilde{a},\tilde{b},\tilde{\beta})
\in G_{\sH\te \sK}\]
by setting   $\tilde{a}=a\te d$ and letting  $\tilde{b}$ 
 be a selfadjoint extension of $a\te d+b\te I$ and with 
 $\tilde{\beta}$ given by certain rather complicated formula 
(see \cite{ax+b}).  Observe, that $(a\te d)\za (b\te I)$
 and that $(a\te d,I\te\de, b\te I,\be\te I)\in M^2_{\sH\te \sK}$. 
It turns out that in order to make the operation $\mlot$ associative 
 one has to assume that $\hb$ is given by \rf{hb} and that 
$k$ describing selfadjoint extension of  $a\te d+b\te I$ (see formula 
\rf{ext} and the explanation below) is related to $\hb$ via \rf{hb}.
 
\section{\Qa\  as a  \csta crossed product}
\label{cp}

We recall that  $A$ denotes the \cstal\ 
of all continuous vanishing at infinity functions on 
\qa  \ . The  \cstal\ $A$ is  generated (see explanation in Section \ref{not})
 by unbounded operators $\log a$, $\;$ $b$ and $i\be b$, 
such that  $(a,b,\be)\in G_{\sH}$.

Assume that $\ker b=\{0\}$. This assumption is not very 
restrictive since every $b$ is a 
direct sum of $b_1$ invertible and $b_2=0$,
 and the case $b_2=0$ is not interesting.  
It was proved in \cite{ax+b} that 
 in that case  the multiplicative unitary operator 
$W\in B(\cH\te \cH)$ related to 
the \qa   is given by
\beq 
W=\fh \left(\eh b^{-1}a \te b,(-1)^k(\be\te \be)
\chi( b\te b <0)\right)^*e^{\frac{i}{\hb}\log (|b|^{-1})\te \log a} .
\label{W}
\eeq
In fact, to have a manageable multiplicative unitary, which is essential 
in S.L. Woronowicz theory of multiplicative unitary 
operators (see \cite{wmul}), one needs a slightly more complicated formula
 for $W$. However,
for our purpose, the formula above is good enough.

From the theory of multiplicative unitaries we know that 
 the structure of the quantum group is encoded in its multiplicative unitary 
operator. More precisely, for any $d\in A$ 
 comultiplication  $\Delta\in \Mor (A,A\te A)$ 
is given by
\[\Delta(d)=W(d\te \id)W^*.\]
Moreover, $\Delta$ may be extended to  unbouned operators 
 affiliated with $A$ and is given on 
  generators  of $A$ by 
\beq
\label{deltaa}
\Delta(a)=W(a\te \id)W^*=a\te a,
\eeq
\beq
\label{deltab}
  \Delta(b)=W(b\te \id)W^*=\left[a\te b+b\te I\right]_{(-1)^k
(\beta\te\beta)\chi(b\otimes b<0)},
\eeq
\[\Delta(i b\be)=W(ib\be\te \id)W^*=i \Delta (b) W(\beta\te\id)W^*\ .\]
By  abuse of notation we also write 
\beq
\label{deltabe}
\Delta(\be)=W(\be\te \id)W^*\ .
\eeq
Thus defined  $\Delta$ is associative

We construct now a  \csta dynamical system. 
Let $\tau\in [0,+\infty[$ and 
\[b(\tau)=\left[\ba{cc}\tau&0\\0&-\tau\ea\right]\mi
\beta (\tau)=\left[\ba{cc}0&\chi (\tau\neq 0)\\\chi (\tau\neq 0)&0\ea\right]
\ .\]
Then  $(b,\beta)\in M_{\C^2}$. 

Let  $M_{2\times 2}(\C)$ denote the set of all   $2\times 2$   
 matrices over  $\C$.
Let $B_o$ be a  an   
algebra  of all  \ci\  functions   $$f: [0,+\infty[\rightarrow 
M_{2\times 2}(\C),$$
 such that 
$\lim_{t\rightarrow\infty}f(t)=0 $
and $f(0)$ is a multiple of unity, i.e.
$$B_o=\{f\in \Ci([0,+\infty[)\te M_{2\times 2}(\C)\;|\;f(0)=zI_{M_{2\times 2}}\,\whe z\in\C\;\;\}\ .$$
Any function $g\in B_o$ is given by
\[(g(b,\beta))(\tau)=(g_1(b)+g_2(b)i\beta)(\tau)
=\left[\ba{cc}g_1(\tau)&i g_2(\tau)\\
-ig_2(-\tau)&g_1(-\tau)\ea\right]\ ,\]
where $g_1,g_2\in \Ci (\R)$, and  $g_2(0)=0$. 
 Then $B_o$  is a   \cstal.
Observe that $b$ and $i\be b$ are affiliated with   $B_o$ and 
$B_o$ is generated by  $b$ and $i\be b$ in the sense explained in Section 
\ref{not}.
One may identify $B_o$  with the \cstal\  of all continuous 
vanishing at infinity  functions
   on  $M_{\C^2}$ (see Section 1.3 of \cite{paper1} and references therein).
Moreover 
$$\M(B_o)=\{f\in C_{\rm bounded} ([0,+\infty[,M_{2\times 2}(\C))\;|\;  
f(0)=zI_{M_{2\times 2}}\ , \whe  z\in\C\;\;\}.$$
Analogously, one may identify $\M(B_o)$  with the \cstal\  of all continuous 
 bounded   functions
   on  $M_{\C^2}$.

Let us introduce an action
$\sigma\in {\rm Aut}(\M (B_o))$\index{$\sigma\;\;$} of   
$\R$ 
given  for any  function
$f\in B_o$ 
by
\[(\sigma_t f)(\tau)=f(e^{\hb t}\tau)\whe \tau\in\R_+ \mi t\in\R\ . \]
Then $(B_o,\R,\sigma)$ is a \csta dynamical system (see e.g. \cite{ped}). 
Let us 
denote the  \csta crossed product algebra
coming with this system by $A_{cp}$
\[A_{cp}=B_o\times_{\sigma} \R \ .\]

 From the definition of  \csta  crossed  product follows that 
$M(A_{cp})$  contains a one parameter, strictly continouous group of unitary 
 operators  implementing  the  action  $\sigma$ of the group $\R$
 on algebra $B_o$. 
Denote the infinitesimal generator of this group by $\log a$. 
Then $a$ is 
a strictly positive operator affilated with $A$.
For any   $t\in \R$, we have that  $a^{it}\in M(A)$
 is a unitary operator  and   
  map  $\R\ni t\rightarrow a^{it}d$ is  
   continouous   for any   $d\in A$.
The action  $\sigma$ of  $\R$ 
 is   given  for any  $f\in B_o$ and $ t\in \R$ by 
\[\sigma_t (f)=a^{it}f a^{-it}\ .\]
It follows that $a\za b$ and $a\be=\be a$.
Therefore
$(a,b,\be)\in G_\cH$.

It is well-known that a linear envelope of the set  
\[ \{\ fg(\log a)\ : f\in B_o, \ g\in\Cir\}\]
 is dense in    $A=B_o\times_\sigma\R$. 
Hence $B_o \subset M(A)$.
It turns out  (\cite[Proposition 3.1]{ax+b}) 
 that  $\log a,\  b,$ and  $\be$ generate  $A_{cp}$.
Moreover,  Proposition 3.2 \cite{ax+b} says that for any triple  
$(\tilde{a},\tilde{b},\tilde{\be})\in G_{\cH}$ there is 
 a unique  representation 
$\pi\in {\rm Rep} (A_{cp},H)$ such that $\tilde{a}=\pi (a)$, 
 $\tilde{b}=\pi(b)$ and  $\tilde{\beta}=\pi (\beta)$. 
It was proven in  \cite{ax+b} that $A_{cp}=A$, i.e. that  
 the crossed  product  algebra $A_{cp}$ is the same as algebra 
 $A$ of all continuous, tending to 0 in infinity functions on 
 the \qa (see \cite{ax+b}). From now  on we will use letter $A$
 for both these algebras.

Now we will  construct a dual action of $\R$ on the crossed 
product algebra related to the dynamical system 
$(B_o,\R,\sigma)$.
To this end let us consider the map 
\beq
\label{teta}
\theta_t =(\id\te\fil _t)\Delta \ ,
\eeq
where the map 
$\fil \in \Mor (A,\Cir)$ 
 is such that \dlad  $t\in\R$
\beq
\fil _t (\log a)=\hb t \mi  \fil _t (b)=0\ .
\label{zgen}
\eeq
From  Proposition 3.2 of  \cite{ax+b} it follows that  there is only one such 
map. 
Observe that   $\theta_t$ is an automorphism of  $A$ and that 
\[\theta_0=\id \mi \theta_s\theta_t=\theta_{s+t}\]
The map $t\rightarrow \theta_t(d)$ is continuous for any  $d\in A$.
Morever, for any  $x,t\in \R$
\[\theta_x(a^{it})=e^{i\hb tx}a^{it}\ . \]
Thus we showed that  $(A,\R,\theta)$ is  
a \csta dynamical system and that it is the system dual to the 
 \csta dynamical system  $(B_o,\R,\sigma)$.

By \rf{deltaa} and  \rf{zgen} we get 
\beq
\label{tetaa}
 (\id\te\fil _t)\Delta (a)=(\id\te\fil _t)(a\te a)=e^{\hb t} a\ .
\eeq
Moreover, by  \rf{W}  and  \rf{zgen}                          
\[(\id\te\fil _t)\fh \left(\eh b^{-1}a \te b,(-1)^k(\be\te \be)
\chi( b\te b <0)\right)^*e^{\frac{i}{\hb}\log (|b|^{-1})\te \log a}
=\fh (0,0)^*e^{it\log (|b|^{-1})}\ .\]
Hence by \rf{deltab}
 \[(\id\te\fil _t)\Delta (b)=b\ .\]
Analogously
 \[(\id\te\fil _t)\Delta (i\be b)=i\be b\ .\]
It follows that for any  $t\in R$ and $g\in \M(B_o)$ we have  
\[\theta_t (g)=g \ .\]
Since $(B_o,\R,\sigma)$ and $(A,\R,\theta)$ are dual
 \csta dynamical systems, it follows that 
 $(B_o^{\prime\prime},\R,\sigma)$ and   
$(A^{\prime\prime},\R,\theta)$ are dual $W^*$-dynamical systems,
 if we only extend  $\sigma$ and  $\theta$ appropiately. 

Let $\cK$ be a Hilbert space, finite- ore infinitedimensional.
Observe that    
$(B(\cK)\te B_o^{\prime\prime},\R,I\te\sigma)$
 is also a  $W^*$-dynamical system  and its crossed  product  
  von Neumann algebra $W^*$ is  
$B(\cK)\te A^{\prime\prime}$.
Analogously, the  $W^*$-dynamical system dual to  
  $(B(\cK)\te B_o^{\prime\prime},\R,I\te\sigma)$ 
 is   $(B(\cK)\te A^{\prime\prime},\R,I\te\theta)$.

In what follows we will need the  Proposition \ref{niezm} below, 
 which is an easy consequence 
of  Theorem  7.10.4 from \cite{ped}:
\bst
\label{niezm}
Let $m\in B(\cK)\te A^{\prime\prime}$ and let for any  
 $t\in\R$  
\[ (\id\te \theta_t)(m)=m\ .\]
Then
\[m\in B(\cK)\te B_o^{\prime\prime}\ .\]
\est
We also need
\bst
Let $B$ be a  \csta -algebra and let $w\in\M(B\te A)$. 
Then the  map   $\R\ni t\rightarrow (\id\te\fil_t) w \in \M(B)$ 
is strictly continuous.
\label{dostone}
\est
\dow
We know that  $\fil \in \Mor(A,\Ci (\R))$ and therefore  $(\id \te \fil)w\in 
\M(B\te \Ci (\R))$. S.L. Woronowicz  in \cite{wunb} showed that  elements of  
$\M(B\te \Ci (\R))$ are bounded, strictly  continuous functions  
on $\R$ with values in  $M(B)$.\hfill\qed

\section{Representation theorem}
\label{ost}
\bde[Unitary representation] 
A unitary operator $V\in \M(CB(\cK\te A)$ is  
 called a (strongly continuous) unitary representation of the \qa \ 
if 
\[W_{23}V_{12}=V_{12}V_{13}W_{23}\ ,\] 
or equivalently
\beq
\label{ru}
(\id\te \Delta)V=V_{12}V_{13}\ .
\eeq
\ede

Observe that in case of  the classical group   condition  
\rf{ru} is equivalent to the 
 classical  definition of   a unitary representation, i.e.  
a representation  $U$   is a map 
$$U:G\ni g\rightarrow U_g\in B(\cK)$$
 such that  $U_g$ is  unitary for any  $g\in G$, and  
 for any $g,h\in G$ we have  $U_g U_h=U_{gh}$. 

It was proven in  \cite{ax+b} that 
\bst
Let  $(a,b,\be)\in G_{\sH}$ and $(c,d,\de)\in G_{\sK}$ and let 
$\ker b=\{0\}$.
Then the  operator $V\in \M(CB(\cH\te A)$ given by 
\beq
\label{V}
V(\log a,b,\be)=\fh \left(d \te b,(\de\te \be)
\chi(d\te b <0)\right)^*e^{\frac{i}{\hb}\log c\te \log a}
\eeq satisfies  
\[(\id \te \Delta)V=V_{12}V_{13}\ ,\]
so  $V$ is a \ru \  of the \qa .
\est
We will prove that all unitary reperesentations of the \qa \ 
are of the form described above. This is the main result of this paper.

\btw
A  $V\in \M(CB(\cH\te A))$ is a unitary representation of 
 the \qa \  on a Hilbert space $\cK$ 
iff there exists  $(c,d,\de)\in G_{\sK}$ such that 
\[V(\log a,b,\be)=\fh \left(d \te b,(\de\te \be)
\chi(d\te b <0)\right)^*e^{\frac{i}{\hb}\log c\te \log a}\ ,\]
 where  $(a,b,\be)\in G_{\sH}$ and $\log a$, $b$ and $ib\be$ 
  are the generators of $A$.
\etw
\dow
Let $V$ be a unitary representation of 
the \qa \  
on a Hilbert space $\cK$. Then for any  $t\in\R $ 
 the operator $(\id\te\fil _t)V\in B(\K)$ is unitary.

Applying  $(\id\te\fil_s\te\fil_t)$ to both sides of 
 \rf{ru} we get 
\[V(\hb (s+t),0,0)=V(\hb s,0,0)V(\hb t,0,0)\]
Hence 
\[(\id\te\fil_{s+t})V=(\id\te\fil_{s})V(\id\te\fil_{t})V\]
i.e.  $(\id\te\fil)V$ is a representation of  $\R$. 

The strict topology   coincides on  $B(H)=M(CB(H))$ 
with the *-strong operator topology.
Since  *-strong operator topology  is stronger than  
strong operator topology, 
 Proposition  \ref{dostone}  shows that the map  
\[ \R\ni t\rightarrow (\id\te\fil_t)V\in B(\K)\]
is strongly  \ci.
 Therefore, by  Stone's Theorem, there is 
 a selfadjoint, strictly positive  operator $c$ acting on  $\K$
 and such that 
\beq
\label{zau1}
(\id\te\fil_t)V=e^{it\log c}
\eeq
\dlad $t\in\R$.

Observe that by \rf{teta}, \rf{ru} and  \rf{zau1}
\[(\id\te\theta_t)V=V(e^{it\log c}\te\id) \ .\]
Moreover, by \rf{tetaa}
\[(\id\te\theta_t)e^{-\frac{i}{\hb}\log c\te\log a}=
(e^{-it\log c}\te\id)e^{-\frac{i}{\hb}\log c\te\log a}\]
 Hence 
\[(\id\te\theta_t)Ve^{-\frac{i}{\hb}\log c\te\log a}=
Ve^{-\frac{i}{\hb}\log c\te\log a}\ .\]
Observe that 
\[Ve^{-\frac{i}{\hb}\log c\te\log a}
\in B(\cK)\te A^{\prime\prime}\ .\] 
By Proposition  \ref{niezm}
\[Ve^{-\frac{i}{\hb}\log c\te\log a}=f(b,\be)\ ,\]
where  $f\in  B(\cK)\te B_o^{\prime\prime}$.
It means that one can  (\cite{wunb}) identify  $f$ 
 with a measurable (operator) function 
 $f:M\rightarrow B(\cK\te\cH)$ . 
Because $V$ is unitary and 
\[V=f(b,\be)e^{\frac{i}{\hb}\log c\te\log a}\ ,\]
it follows that  $f$ may be identified 
 with  a \ci\ and  bounded (operator) function  with values in  
unitary operators acting on  $\cK$. 
Let us compute 
\[(\id\te\Delta)V=(\id\te\Delta)\left(f(b,\be)e^{\frac{i}{\hb}
\log c\te\log a}\right)=f(\Delta(b),\Delta(\be))e^{\frac{i}{\hb}
\log c\te\log \Delta( a)}=\]
\beq
\label{1}
=f([a\te b+b\te I]_{\tau},\Delta(\be))e^{\frac{i}{\hb}
\log c\te\log ( a\te a)}\ .
\eeq
Applying $(\id \te \fil_t\te \id)$ to both sides of \rf{ru}
 we get 
\[(\id \te \fil_t\te \id)V_{12}V_{13}=(e^{it\log c}\te\id)V\]
and
\[(\id \te \fil_t\te \id)(\id\te\Delta)V=f(e^{\hb t}b,\be)
(e^{\frac{i}{\hb}\log c\te\log a})
(e^{it\log c}\te\id)\ .\]
Comparing the resulting expressions we obtain
\beq
\label{por1}
(e^{it\log c}\te\id)V
=f(e^{\hb t}b,\be)
(e^{\frac{i}{\hb}\log c\te\log a})
(e^{it\log c}\te\id)\ .
\eeq
On the other hand, applying   $(\id \te\id \te \fil_t)$ to both 
sides of \rf{ru}
 we get 
\[(\id \te \id \te \fil_t)V_{12}V_{13}=V(e^{it\log c}\te\id)\]
and 
\[(\id \te \id \te \fil_t)(\id\te\Delta)V=f(b,\be)
(e^{\frac{i}{\hb}\log c\te\log a})
(e^{it\log c}\te\id).\]
Comparing the resulting expressions we obtain
\beq
\label{por2}
V(e^{it\log c}\te\id)=f(b,\be)
(e^{\frac{i}{\hb}\log c\te\log a})
(e^{it\log c}\te\id)\ .
\eeq
We insert now  $\log a$ instead of  $\hb t$ 
 in formulas \rf{por1} and \rf{por2} we get 
\[(e^{\frac{i}{\hb}
\log c \te \log a}\te\id)V_{13}
=f(a\te b,I\te \be)
(e^{\frac{i}{\hb}\log c\te\log a})_{13}
(e^{\frac{i}{\hb}
\log c \te \log a})_{12}\]
\[V_{12}(e^{\frac{i}{\hb}
\log c \te \log a})_{13}=f(b\te I,\be\te I)
(e^{\frac{i}{\hb}\log c\te\log a})_{12}
(e^{\frac{i}{\hb}
\log c \te \log a})_{13}\ .\]
Then
\beq
\label{2}
V_{12}=f(b\te I,\be\te I)(e^{\frac{i}{\hb}
\log c \te a})_{12}
\eeq
\beq
V_{13}= (e^{\frac{i}{\hb}
\log c \te a})_{12}^* f(a\te b,I\te \be)(e^{\frac{i}{\hb}
\log c \te a})_{13}(e^{\frac{i}{\hb}
\log c \te a})_{12}\ .
\label{565}
\eeq
remembering that $V$ satisfies  $(\id \te\Delta)V= V_{12}V_{13}$ 
and using \rf{1}, \rf{2} and \rf{565}, we get
\[
f([a\te b+b\te I]_{\tau},\Delta(\be))e^{\frac{i}{\hb}
\log c\te\log ( a\te a)}=f(b\te I,\be\te I) f(a\te b,I\te \be)(e^{\frac{i}{\hb}
\log c \te \log a})_{13}(e^{\frac{i}{\hb}
\log c \te \log a})_{12}\ .
\]
Moreover
\[e^{\frac{i}{\hb}
\log c\te\log ( a\te a)}=e^{\frac{i}{\hb}
\log c\te\{\log ( I\te a)+\log ( a\te I)\}}=
(e^{\frac{i}{\hb}
\log c \te \log a})_{13}(e^{\frac{i}{\hb}
\log c \te \log a})_{12}\ .\]
Therefore
\[f([a\te b+b\te I]_{\tau},\Delta(\be))=
f(b\te I,\be\te I) f(a\te b,I\te \be)\]
or equivalently
\beq
\label{rost}
f([a\te b+b\te I]_{\tau},\Delta(\be))^*=
f(a\te b,I\te \be)^*f(b\te I,\be\te I)^*\ .
\eeq
Using  \rf{deltabe} and  \rf{W} and remembering that    
$\be$ commutes  with $s |b|^{-1}$, we get
\beq\Delta (\be)=
\label{deltabe2}
\eeq
\[\!\!\!\!\!\!\!\!\!=\fh \left(\eh b^{-1}a \te b,(-1)^k(\be\te \be)
\chi( b\te b <0)\right)^*
(\be\te I)
\fh \left(\eh b^{-1}a \te b,(-1)^k(\be\te \be)
\chi( b\te b <0)\right)\ .\]
Let us introduce the notation 
\beq
\label{oznrs}
R=a \te b\ , \od \rho=I\te \be \od 
 S= b\te I\ ,\od \sigma=\be\te I\ .
\eeq
We see that   $(R,\rho,S,\sigma)\in M^2_{\cH\te \cH}$.
Moreover, if we define  
\beq
\label{ozntt}
T=\eh S^{-1}R=\eh b^{-1}a\te b \mi
\tau = (-1)^k(\be \te \be)
\ ,
\eeq
then also  $(T,\tau)\in M_{\cH}$.
Since $a\za b$, then  by  Theorem 6.1 from  \cite{qexp} 
\[\sign ( \eh b^{-1}a )=\sign (b)\; \sign (a)\ .\]
 Since  $a>0$, we finally get 
\[\chi(T<0)=\chi( \eh b^{-1}a\te b <0)=\chi(  b\te b <0)\ .
\label{znaki}
\]
Inserting \rf{oznrs} 
 and  \rf{ozntt} into formula  \rf{deltabe2}
 and using  \rf{znaki}  we get 
\[\Delta (\be)=\fh(T,\tau \chi(T<0) )^*\sigma \fh(T,\tau \chi(T<0) )\ . \]
Moreover, we observe that the following selfadjoint extensions 
are the same
\[
[a\te b+b\te I]_{\tau}=
[R+S]_{\tau\chi(T<0)}
\ .\]
Using two last formulas and substituting  \rf{oznrs} into the right 
hand side of \rf{rost}, one can write \rf{rost} 
\[
f((R,\rho) \mlot _M  (S,\sigma))^*
=f(R,\rho)^*f(S,\sigma)^*.\]
From Theorem 3.3  \cite{paper1} (recalled in Section 2) follows that if 
  a function $f$ is a measurable (operator) function 
 $f:M\rightarrow B(\cK\te\cH)$  and satisfies the above condition, 
 then  it is given by 
\[ f(b,\be)=\fh \left( d\te b,(\delta\te \beta)\chi(d\te b <0) 
\right)^* \ ,\]
where $( d,\delta)\in M\sK$.

Therefore
\beq
\label{hura}
V=\fh \left( d\te b,(\delta\te \beta)\chi(d\te b <0) 
\right)^*e^{\frac{i}{\hb}
\log c \te \log a}
\eeq
What is left to prove is that 
\[c \za d \mi \de c= c \de \ .\]
To this end, let us observe that combining \rf{hura} and \rf{565}
 we get
\[\fh (d\te I\te b,(\delta\te I\te\be)\chi
(d\te I\te b<0))^*= \]
\[=(e^{\frac{i}{\hb}
\log c \te \log a})_{12}^* \fh(d\te a\te b,(\delta\te I\te \be)\chi
(d\te I\te b<0))^*(e^{\frac{i}{\hb}
\log c \te \log a})_{12}\ .\]
Hence
\beq 
\label{stad1} (d\te \id)=
e^{-\frac{i}{\hb}\log c\te\log a }
(d\te a) 
e^{\frac{i}{\hb}\log c\te\log a }
\eeq
and
\beq
 \label{stad2} (\delta\te \id)=
e^{-\frac{i}{\hb}\log c\te\log a }
(\delta\te I) 
e^{\frac{i}{\hb}\log c\te\log a }\ .
\eeq
Substituting   $a=e^{\hb k}I$ in  formula  \rf{stad1} 
 we get  $c\za d$.
Moreover, after the same substitution in formula \rf{stad2} we see 
 that  
  $c$ commutes with  $\delta$.
Hence  $(c,d,\delta)\in G\sK$, which completes the proof.\hfill\qed

\section*{Acknowledgments}

This is a part of the author's Ph.D. thesis \cite{PhD} written under 
the supervision of Professor Stanis\l aw L. Woronowicz at the 
Department of mathematical Methods in Physics at the Warsaw University.
 The author is greatly indebted to Professor Stanis\l aw L. Woronowicz 
for stimulating discussions and important hints and comments. The author also
 wishes to thank Professor Wies\l aw Pusz and Marek Bo\.zejko for several 
helpful suggestions and Piotr So\l tan for reading carefully the manuscript.


\end{document}